\input ./style/arxiv-general.cfg
\documentclass[aos,MSNbibl,nameyear,seceqn,dvips]{arximspdf}
\makeatletter
   \@ifpackageloaded{graphicx}{}{\usepackage{graphicx}}
\makeatother

\doi{10.1214/15-AOS1375}
\volume{44}
\issue{5}
\pubyear{2016}
\firstpage{1931}
\lastpage{1956}
\docsubty{FLA}

\makeatletter
\newcommand{\eqref}[1]{(\ref{#1})}

\newtheorem{theorem}{Theorem}[section]
\newproclaim{assumption}{Assumption}[section]

\newproclaim{definition}{Definition}[section]
\newproclaim{example}{Example}[section]
\newproclaim{remark}{Remark}[section]

\newcommand{\de}{\delta}
\newcommand{\e}{\mathbb{E}}
\newcommand{\p}{\mathbb{P}}
\newcommand{\var}{\operatorname{var}}
\newcommand{\la}{\lambda}
\newcommand{\wh}{\widehat}
\newcommand{\mX}{\mathbf{X}}
\newcommand{\mY}{\mathbf{Y}}

\def\bX{\mathbf{X}}
\def\bY{\mathbf{Y}}
\def\bzero{\mathbf{0}}

\newcommand{\cov}{\operatorname{cov}}
\newcommand{\corr}{\operatorname{corr}}

\makeatother

\begin{document}
\begin{frontmatter}

\title{Cram\'er-type moderate deviations for Studentized two-sample
$U$-statistics with applications}
\runtitle{Studentized two-sample $U$-statistics}

\begin{aug}
\author[A]{\fnms{Jinyuan} \snm{Chang}\thanksref{t1,m4, m1}\ead[label=e1]{jingyuan.chang@unimelb.edu.au}},
\author[B]{\fnms{Qi-Man} \snm{Shao}\thanksref{t2,m2}\ead[label=e2]{qmshao@cuhk.edu.hk}}
\and
\author[C]{\fnms{Wen-Xin} \snm{Zhou}\corref{}\thanksref{t3,m3,m1}\ead
[label=e3]{wenxinz@princeton.edu}}
\runauthor{J. Chang, Q.-M. Shao and W.-X. Zhou}

\affiliation{Southwestern University of Finance and Economics,\thanksmark{m4} University of Melbourne,\thanksmark{m1} Chinese University
of Hong Kong\thanksmark{m2} and
Princeton University\thanksmark{m3}}

\thankstext{t1}{Supported in part by the Fundamental Research Funds for the Central Universities (Grant No. JBK160159,
JBK150501), NSFC (Grant No. 11501462), the Center of Statistical Research at SWUFE and the Australian
Research Council.}
\thankstext{t2}{Supported by Hong Kong Research Grants Council GRF 603710 and 403513.}
\thankstext{t3}{Supported by NIH Grant R01-GM100474-4 and a grant
from the Australian Research Council.}

\footnotetext{\textbf{Tribute}: Peter was a brilliant and prolific researcher,
who has made enormously influential contributions to mathematical statistics and probability theory.
Peter had extraordinary knowledge of analytic techniques that he often applied with ingenious
simplicity to tackle complex statistical problems. His work and service have had a profound
impact on statistics and the statistical community. Peter was a generous mentor and friend with
a warm heart and keen to help the young generation. Jinyuan Chang and Wen-Xin Zhou are extremely
grateful for the opportunity to learn from and work with Peter in the last two years at the
 University of Melbourne. Even in his final year, he had afforded time to guide us.
 We will always treasure the time we spent with him. Qi-Man Shao is so grateful
for all the helps and supports that Peter had provided during the various stages of his career.
Peter will be dearly missed and forever remembered as our mentor and friend.}

\address[A]{J. Chang\\
School of Statistics\\
Southwestern University of Finance\\
\quad and Economics\\
Chengdu, Sichuan 611130\\
China \\
and\\
School of Mathematics and Statistics\\
University of Melbourne\\
Parkville, Victoria 3010\\
Australia \\
\printead{e1}}
\address[B]{Q.-M. Shao\\
Department of Statistics\\
Chinese University of Hong Kong \\
Shatin, NT \\
Hong Kong \\
\printead{e2}}
\address[C]{W.-X. Zhou\\
Department of Operations Research\\
\quad and Financial Engineering \\
Princeton University \\
Princeton, New Jersey 08544 \\
USA\\
and\\
School of Mathematics and Statistics \\
University of Melbourne \\
Parkville, Victoria 3010 \\
Australia\\
\printead{e3}}
\end{aug}

%
\received{\smonth{6} \syear{2015}}


\begin{abstract}
Two-sample $U$-statistics are widely used in a broad range of
applications, including those in the fields of biostatistics and
econometrics. In this paper, we establish sharp Cram\'er-type moderate
deviation theorems for Studentized two-sample $U$-statistics in a
general framework, including the two-sample $t$-statistic and
Studentized Mann--Whitney test statistic as prototypical examples. In
particular, a refined moderate deviation theorem with second-order
accuracy is established for the two-sample $t$-statistic. These results
extend the applicability of the existing statistical methodologies from
the one-sample $t$-statistic to more general nonlinear statistics.
Applications to two-sample large-scale multiple testing problems with
false discovery rate control and the regularized bootstrap method are
also discussed.\vspace*{-3pt}
\end{abstract}

%
\end{frontmatter}

\section{Introduction} \label{intro.sec}

The $U$-statistic is one of the most commonly used nonlinear and
nonparametric statistics, and its asymptotic theory has been well
studied since the seminal paper of \citet{Hoeffding1948}. $U$-statistics
extend the scope of parametric estimation to more complex nonparametric
problems and provide a general theoretical framework for statistical
inference. We refer to \citet{KoroljukBorovskich1994} for a systematic
presentation of the theory of $U$-statistics, and to \citet
{KowalskiTu2007} for more recently discovered methods and contemporary
applications of $U$-statistics.

Applications of $U$-statistics can also be found in high dimensional
statistical inference and estimation, including the simultaneous
testing of many different hypotheses, feature selection and ranking,
the estimation of high dimensional graphical models and sparse, high
dimensional signal detection. In the context of high dimensional
hypothesis testing, for example, several new methods based on
\mbox{$U$-statistics} have been proposed and studied in \citet{ChenQin2010},
\citet{ChenZhangZhong2010} and \citet{ZhongChen2011}. Moreover, \citet
{LiPengZhangZhu2012} and \citeauthor{LiZhongZhu2012}
(\citeyear{LiZhongZhu2012}) employed $U$-statistics
to construct independence feature screening procedures for analyzing
ultrahigh dimensional data.

Due to heteroscedasticity, the measurements across disparate subjects
may differ significantly in scale for each feature. To standardize for
scale, unknown nuisance parameters are always involved and a natural
approach is to use Studentized, or self-normalized statistics. The
noteworthy advantage of Studentization is that compared to standardized
statistics, Studentized ratios take heteroscedasticity into account and
are more robust against heavy-tailed data. The theoretical and
numerical studies in \citet{DelaigleHallJin2011} and \citeauthor{ChangTangWu2013} (\citeyear
{ChangTangWu2013,ChangTangWu2014}) evidence the importance of using
Studentized statistics in high dimensional data analysis. As noted in
\citet{DelaigleHallJin2011}, a careful study of the moderate deviations
in the Studentized ratios is indispensable to understanding the common
statistical procedures used in analyzing high dimensional data.

Further, it is now known that the theory of Cram\'er-type moderate
deviations for Studentized statistics quantifies the accuracy of the
estimated $p$-values, which is crucial in the study of large-scale
multiple tests for controlling the false discovery rate [\citet
{FanHallYao2007,LiuShao2010}]. In particular, Cram\'er-type moderate
deviation results can be used to investigate the robustness and
accuracy properties of $p$-values and critical values in multiple
testing procedures. However, thus far, most applications have been
confined to $t$-statistics [\citet
{FanHallYao2007,WangHall2009,DelaigleHallJin2011,CaoKosorok2011}]. It is
conjectured in \citet{FanHallYao2007} that analogues of the theoretical
properties of these statistical methodologies remain valid for other
resampling methods based on Studentized statistics. Motivated by the
above applications, we are attempting to develop a unified theory on
moderate deviations for more general Studentized nonlinear statistics,
in particular, for two-sample $U$-statistics.

The asymptotic properties of the standardized $U$-statistics are
extensively studied in the literature, whereas significant developments
are achieved in the past decade for one-sample Studentized
$U$-statistics. We refer to \citet{WangJingZhao2000} and the references
therein for Berry--Esseen-type bounds and Edgeworth expansions. The
results for moderate deviations can be found in \citet
{VandemaeleVeraverbeke1985}, \citet{LaiShaoWang2011} and \citet
{ShaoZhou2012}. The results in \citet{ShaoZhou2012} paved the way for
further applications of statistical methodologies using Studentized
$U$-statistics in high dimensional data analysis.

Two-sample $U$-statistics are also commonly used to compare the
different (treatment) effects of two groups, such as an experimental
group and a control group, in scientifically controlled experiments.
However, due to the structural complexities, the theoretical properties
of the two-sample $U$-statistics have not been well studied. In this
paper, we establish a Cram\'er-type moderate deviation theorem in a
general framework for Studentized two-sample $U$-statistics, especially
the two-sample $t$-statistic and the Studentized Mann--Whitney test. In
particular, a refined moderate deviation theorem with second-order
accuracy is established for the two-sample $t$-statistic.

The paper is organized as follows. In Section~\ref{cmd.U}, we present
the main results on Cram\'er-type moderate deviations for Studentized
two-sample $U$-statistics as well as a refined result for the
two-sample $t$-statistic. In Section~\ref{HT.sec}, we investigate
statistical applications of our theoretical results to the problem of
simultaneously testing many different hypotheses, based particularly on
the two-sample $t$-statistics and Studentized Mann--Whitney tests.
Section~\ref{sim.sec} shows numerical studies. A discussion is given in
Section~\ref{se:discu}. All the proofs are relegated to the
supplementary material [\citet{supp}].

\section{Moderate deviations for Studentized $U$-statistics}
\label{cmd.U}

We use the following notation throughout this paper. For two sequences
of real numbers $a_n$ and $b_n$, we write $a_n \asymp b_n$ if there
exist two positive constants $c_1$, $c_2$ such that $c_1 \leq a_n/ b_n
\leq c_2$ for all $n\geq1$, we write $a_n= O(b_n)$ if there is a
constant $C$ such that $|a_n|\leq C |b_n|$ holds for all sufficiently
large $n$, and we write $a_n \sim b_n$ and $a_n=o(b_n)$, respectively,
if $\lim_{n\rightarrow\infty} a_n/b_n =1$ and $\lim_{n\rightarrow
\infty} a_n/b_n =0$. Moreover, for two real numbers $a$ and $b$, we
write for ease of presentation that $a \vee b=\max(a, b)$ and $a
\wedge
b=\min(a, b)$.

\subsection{A review of Studentized one-sample $U$-statistics}
\label{One.U}

We start with a brief review of Cram\'er-type moderate deviation for
Studentized one-sample $U$-statistics. For an integer $s \geq2$ and
for $n> 2s$, let $X_1,\ldots, X_n$ be independent and identically
distributed (i.i.d.) random variables taking values in a metric space
$(\mathbb{X}, \mathcal G)$, and let $h: \mathbb{X}^d \mapsto\mathbb
{R}$ be a symmetric Borel measurable function. Hoeffding's
$U$-statistic with a kernel $h$ of degree $s$ is defined as
\[
U_n = \frac{1}{{n \choose s}} \sum_{1\leq i_1 < \cdots< i_s \leq n}
h(X_{i_1},\ldots, X_{i_s}), \label{u-stat}
\]
which is an unbiased estimate of $\theta=\e\{ h(X_1,\ldots,X_s)\}$.
In particular, we focus on the case where $\mathbb{X}$ is the Euclidean
space $\mathbb{R}^r$ for some integer $r\geq1$. When $r\geq2$, write
$X_i =(X_i^1,\ldots, X_i^r)^{\mathrm{T}}$ for $i=1,\ldots, n$.

Let
\[
h_1(x)=\e\bigl\{ h(X_1,\ldots,X_s) |
X_1=x \bigr\} \qquad \mbox{for any } x =\bigl(x^1, \ldots,
x^r\bigr)^{\mathrm{T}}\in\mathbb{R}^r
\]
and
\[
\sigma^2 =\var\bigl\{ h_1(X_1) \bigr\},\qquad
  v_h^2 = \var\bigl\{h(X_1,
X_2,\ldots ,X_s) \bigr\}. \label{var.def}
\]
Assume that $0< \sigma^2 < \infty$, then the standardized
nondegenerate $U$-statistic is given by
\[
Z_n = {{n}^{1/2} \over s \sigma} (U_n-\theta).
\label{standU}
\]


Because $\sigma$ is usually unknown, we are interested in the following
Studentized $U$-statistic:
%
\begin{equation}
\wh{U}_n ={ {n}^{1/2} \over s \wh{\sigma
} } (U_n-\theta),
\label{stud-u}
\end{equation}
where $ \wh{\sigma}^2$ denotes the
leave-one-out jackknife estimator of $\sigma^2$ given by
\begin{eqnarray*}
\wh{\sigma}^2 &=& \frac{(n-1)}{(n-s)^2} \sum
_{i=1}^n(q_i-U_n)^2
\qquad\mbox{with} \label{s1}
\\
q_i &= & \frac{1}{{n-1 \choose s - 1}} \mathop{\sum
_{1\leq\ell
_1<\cdots< \ell_{s-1} \leq n }}_{
\ell_j\neq i\ \mathrm{for\ each}\ j=1,\ldots, s-1}h(X_i,X_{\ell
_1},
\ldots ,X_{\ell_{s-1}}).
\end{eqnarray*}

\citet{ShaoZhou2012} established a general Cram\'er-type moderate
deviation theorem for Studentized nonlinear statistics, in particular
for Studentized $U$-statistics.

\begin{theorem} \label{thmm.stu.U}
Assume that $ v_p: = [ \e\{| h_1(X_1)-\theta|^p\} ]^{1/p} < \infty$
for some $2 < p \leq3$. Suppose that there are constants $c_0\geq1$
and $\kappa\geq0$ such that for all $x_1,\ldots, x_s \in\mathbb{R}$,
%
\begin{equation}
\bigl\{ h(x_1,\ldots, x_s) - \theta\bigr\}^2
\leq c_0 \Biggl[ \kappa\sigma^2 + \sum
_{i=1}^s \bigl\{ h_1(x_i)-
\theta\bigr\}^2 \Biggr]. \label{k-c}
\end{equation}
Then there exist constants $C, c >0$ depending only on $d$ such that
\[
\frac{ \p( \wh{U}_n\geq x)}{1-\Phi(x)}=1 + O(1) \bigl\{ (v_p/ \sigma)^p (1
+x)^{p} n^{ 1- p/2} +\bigl(a_s^{1/2}+v_h/
\sigma\bigr) (1 +x)^3 n^{-1/2} \bigr\} \label{t1.1a}
\]
holds uniformly for $ 0 \leq x \leq c \min\{ (\sigma/v_p)
n^{1/2-1/p}, (n/a_s)^{1/6} \}$, where $|O(1)|\leq C$ and $a_s=\max(c_0
\kappa, c_0+s )$. In particular, we have
\[
{ \p(\wh{U}_n\geq x) \over1-\Phi(x)
} \to1 \label{t1.1b}
\]
holds uniformly in $x\in[0, o(n^{1/2-1/p}) )$.
\end{theorem}

Condition \eqref{k-c} is satisfied for a large class of $U$-statistics.
Below are some examples.\vspace*{6pt}

\noindent\hspace*{-2pt}{\fontsize{9}{11}\selectfont{
\tabcolsep=0pt
\begin{tabular*}{\textwidth}{@{\extracolsep{\fill}}lcc@{}c@{}}
\hline
\textbf{Statistic} & \textbf{Kernel function} & $\bolds{c_0}$ & $\bolds{\kappa}$ \\
\hline
$t$-statistic & $h(x_1, x_2)=0.5(x_1+x_2)$ & \phantom{0}2 & 0\\
Sample variance & $h(x_1,x_2)=0.5(x_1-x_2)^2$ & 10 & $ (\theta/\sigma
)^2$\\
Gini's mean difference & $h(x_1,x_2)=|x_1-x_2|$ & \phantom{0}8 & $ (\theta/\sigma
)^2$ \\
One-sample Wilcoxon's statistic & $h(x_1,x_2)= I\{x_1+x_2\leq0\}$ & \phantom{0}1
&$\sigma^{-2}$ \\
Kendall's $\tau$ & $h({x}_1, {x}_2)= 2 I\{
(x_2^2-x_1^2)(x_2^1-x_1^1)>0\} $ & \phantom{0}1 & $\sigma^{-2}$ \\
\hline
\end{tabular*}}}

\subsection{Studentized two-sample $U$-statistics}

Let $\mathcal{X}= \{ X_1,\ldots, X_{n_1} \}$ and $\mathcal{Y} =  \{
Y_1, \ldots, Y_{n_2} \}$ be two independent random samples, where
$\mathcal{X}$ is drawn from a probability distribution $P$ and
$\mathcal
{Y}$ is drawn from another probability distribution $Q$. With $s_1$ and
$s_2$ being two positive integers, let
\[
h(x_1,\ldots, x_{s_1}; y_1,\ldots,
y_{s_2})
\]
be a kernel function of order $(s_1, s_2)$ which is real and symmetric
both in its first $s_1$ variates and in its last $s_2$ variates. It is
known that a nonsymmetric kernel can always be replaced with a
symmetrized version by averaging across all possible rearrangements of
the indices.

Set $ \theta:= \e\{ h(X_1,\ldots, X_{s_1}; Y_1,\ldots, Y_{s_2})\}$,
and let
\[
U_{ \bar{n} } = \frac{1}{{n_1 \choose s_1}{n_2 \choose s_2}}\sum_{
1\leq i_1<\cdots<i_{s_1}\leq n_1 } \sum
_{ 1\leq j_1<\cdots
<j_{s_2}\leq
n_2 } h(X_{i_1},\ldots, X_{i_{s_1}};
Y_{j_1},\ldots, Y_{j_{s_2}}), \label{two.U.gen}
\]
be the two-sample $U$-statistic, where $\bar{n} =(n_1, n_2)$. To
lighten the notation, we write $\mX_{i_1, \ldots, i_{\ell}}=(X_{i_1},
\ldots, X_{i_{\ell}})$, $\mY_{j_1, \ldots, j_{k}}=(Y_{j_1}, \ldots,
Y_{j_{k}})$ such that
\[
h(\mX_{i_1, \ldots, i_{\ell}}; \mY_{j_1, \ldots,
j_{k}})=h(X_{i_1},\ldots,
X_{i_{\ell}}; Y_{j_1},\ldots, Y_{j_{k}}),
\]
and define
%
\begin{eqnarray}\label{h1h2}
h_{1}(x)& =& \e\bigl\{ h(\mX_{1,\ldots, s_1} ; \mY_{1, \ldots, s_2}) |
X_{1} =x \bigr\} ,
\nonumber
\\[-8pt]
\\[-8pt]
\nonumber
h_{2}(y) &=& \e\bigl\{ h(\mX_{1,\ldots, s_1} ; \mY_{1, \ldots, s_2}) |
Y_{1}=y \bigr\}.
\end{eqnarray}
Also let $v_h^2 = \var\{ h(\mX_{1,\ldots, s_1} ; \mY_{1, \ldots,
s_2})\} $, $\sigma_{1 }^2=\var\{h_1(X_i)\}$, $\sigma_{2 }^2=\var\{
h_2(Y_j)\}$ and
%
\begin{equation}
\sigma^2 = \sigma_{1 }^2 +
\sigma_{2 }^2 ,\qquad  \sigma_{\bar{n}}^2=
s_1^2 \sigma_1^2
n_1^{-1} +s_2^2
\sigma_2^2 n_2^{-1}. \label{sigma}
\end{equation}

For the standardized two-sample $U$-statistic of the form $\sigma
_{\bar{n}
}^{-1}(U_{\bar{n}} -\theta) $, a uniform Berry--Esseen bound of order
$O\{
(n_1 \wedge n_2)^{-1/2} \}$ was obtained by \citet{HelmersJanssen1982}
and \citet{Borovskich1983}. Using a concentration inequality approach,
\citet{ChenShao2007} proved a refined uniform bound and also established
an optimal nonuniform Berry--Esseen bound. For large deviation
asymptotics of two-sample $U$-statistics, we refer to \citet
{NikitinPonikarov2006} and the references therein.

Here, we are interested in the following Studentized two-sample $U$-statistic:
%
\begin{equation}
\wh{U}_{ \bar{n}} = \widehat{\sigma}_{\bar{n}}^{-1} (
U_{ \bar
{n}} -\theta)  \qquad \mbox{with } \widehat{\sigma}^2_{\bar{n}}
= s_1^2 \widehat{\sigma }_{1}^2
n_1^{-1} + s_2^2 \widehat{
\sigma}_{2}^2 n_2^{-1} ,
\label{two.stu.U}
\end{equation}
where
\[
\wh{\sigma}_{1}^2 = \frac{1}{n_1-1} \sum
_{i=1}^{n_1} \Biggl( q_i -
\frac
{1}{n_1} \sum_{i=1}^{n_1}
q_i \Biggr)^2,  \qquad \wh{\sigma}_{2}^2
= \frac
{1}{n_2-1} \sum_{j=1}^{n_2} \Biggl(
p_j - \frac{1}{n_2} \sum_{j=1}^{n_2}
p_j \Biggr)^2
\]
and
\begin{eqnarray*}
q_i &= &\frac{1}{{n_1-1 \choose s_1-1}{n_2 \choose s_2}}
\mathop{\sum_{ 1\leq i_2 < \cdots< i_{s_1} \leq n_1 }}_{ i_{\ell} \neq i,
\ell=2, \ldots, s_1 }
\sum_{ 1\leq j_1 < \cdots< j_{s_2} \leq n_2 } h (\mX_{ i, i_2,\ldots, i_{s_1} }; \mY_{j_1, \ldots, j_{s_2}}
),
\\
p_j &=& \frac{1}{{n_1 \choose s_1} {n_2-1 \choose s_2-1}} \sum_{1\leq i_1 < \cdots< i_{s_1} \leq n_1 }
\mathop{\sum_{ 1\leq j_2<\cdots<
j_{s_2} \leq n_2 }}_{ j_k \neq j, k=2,\ldots, s_2 } h (\mX_{i_1,
\ldots, i_{s_1}} ;
\mY_{ j, j_2, \ldots, j_{s_2} } ).
\end{eqnarray*}
Note that $\wh{\sigma}_{1}^2$ and $\wh{\sigma}_{2}^2$ are leave-one-out
jackknife estimators of $\sigma_{1}^2$ and $\sigma_{2}^2$, respectively.

\subsubsection{Moderate deviations for \texorpdfstring{$\wh{U}_{\bar{n}}$}{$widehat{U}_{bar{n}}$}}

For $p>2$, let
%
\begin{equation}
v_{1,p} = \bigl[ \e\bigl\{\bigl| h_1(X_1)-
\theta\bigr|^p \bigr\}\bigr]^{1/p} \quad\mbox{and}\quad v_{2,p}
= \bigl[ \e\bigl\{\bigl| h_2(Y_1)-\theta\bigr|^p \bigr
\}\bigr]^{1/p}. \label{p.moment}
\end{equation}
Moreover, put
\[
s = s_1\vee s_2,\qquad \bar{n}=(n_1,
n_2), \qquad n =n_1\wedge n_2 \label{size.not}
\]
and
\[
\la_{\bar{n}} =v_h \biggl( \frac{n_1+n_2}{\sigma_{1 }^2 n_2 +
\sigma_{2 }^2
n_1}
\biggr)^{1/2} \qquad\mbox{with } v_h^2 =\var\bigl\{h(
\mX_{1, \ldots,
s_1} ; \mY_{1,\ldots
, s_2})\bigr\}. \label{la.def}
\]
The following result gives a Cram\'er-type moderate deviation for
$\widehat{U}_{\bar{n}}$ given in \eqref{two.stu.U} under mild assumptions.
A self-contained proof can be found in the supplementary material [\citet{supp}].

\begin{theorem} \label{thmm.1}
Assume that there are constants $c_0 \geq1$ and $\kappa\geq0$ such that
%
\begin{equation}
\bigl\{ h({\mathbf x} ; {\mathbf y} ) -\theta\bigr\}^2 \leq
c_0 \Biggl[ \kappa \sigma^2 + \sum
_{i=1}^{s_1} \bigl\{ h_{1}(x_i)-
\theta\bigr\}^2 + \sum_{j=1}^{s_2}
\bigl\{ h_{2}(y_j)-\theta\bigr\}^2 \Biggr]
\label{ker.reg}
\end{equation}
for all ${\mathbf x}=(x_1, \ldots, x_{s_1})$ and ${\mathbf y} = (y_1,
\ldots,
y_{s_2})$, where $\sigma^2$ is given in \eqref{sigma}. Assume that
$v_{1,p}$ and $v_{2,p}$ are finite for some $2< p \leq3$. Then there
exist constants $C, c>0$ independent of $n_1$ and $n_2$ such that
%
\begin{eqnarray}\label{result.1}
&& \frac{\p(\wh{U}_{\bar{n}} \geq x)}{1- \Phi(x)}
\nonumber
\\[-8pt]
\\[-8pt]
\nonumber
& &\qquad= 1+O(1) \Biggl\{ \sum_{\ell=1}^2
\frac{ v_{\ell,p}^p (1 + x)^p }{
\sigma
_{\ell}^p   n_{\ell}^{ p/2 -1} } + \bigl( a_d^{1/2} + \la_{\bar{n}}
\bigr) (1 + x)^3 \biggl(\frac{n_1+n_2}{n_1n_2} \biggr)^{1/2}
\Biggr\}
\nonumber
\end{eqnarray}
holds uniformly for
\[
0 \leq x \leq c \min \bigl[ (\sigma_{1}/ v_{1,p})
n_1^{ p/2-1} , (\sigma _{2} / v_{2,p})
n_2^{ p/2-1}, a_s^{-1/6} \bigl
\{n_1 n_2/(n_1+n_2)\bigr
\}^{1/6} \bigr], \label{range.x}
\]
where $|O(1)|\leq C$ and $a_s =\max(c_0 \kappa, c_0+s )$. In
particular, as $n \rightarrow\infty$,
%
\begin{equation}
\frac{\p(\wh{U}_{\bar{n}} \geq x )}{1- \Phi(x)} \rightarrow1 \label{cmd}
\end{equation}
holds uniformly in $x \in[0, o(n^{1/2-1/p}) )$.
\end{theorem}

Theorem~\ref{thmm.1} exhibits the dependence between the range of
uniform convergence of the relative error in the central limit theorem
and the optimal moment conditions. In particular, if $p=3$, the region
becomes $0\leq x\leq O( n^{1/6})$. See Theorem~2.3 in \citet
{JingShaoWang2003} for similar results on self-normalized sums. Under
higher order moment conditions, it is not clear if our technique can be
adapted to provide a better approximation for the tail probability $\p(
\wh{U}_{\bar{n}} \geq x )$ for $x$ lying between $n^{1/6}$ and $n^{1/2}$
in order.

It is also worth noticing that many commonly used kernels in
nonparametric statistics turn out to be linear combinations of the
indicator functions and, therefore, satisfy condition \eqref{ker.reg}
immediately.

\subsubsection{Two-sample $t$-statistic}

As a prototypical example of two-sample $U$-statistics, the two-sample
$t$-statistic is of significant interest due to its wide applicability.
The advantage of using $t$-tests, either one-sample or two-sample, is
their high degree of robustness against heavy-tailed data in which the
sampling distribution has only a finite third or fourth moment. The
robustness of the $t$-statistic is useful in high dimensional data
analysis under the sparsity assumption on the signal of interest. When
dealing with two experimental groups, which are typically independent,
in scientifically controlled experiments, the two-sample $t$-statistic
is one of the most commonly used statistics for hypothesis testing and
constructing confidence intervals for the difference between the means
of the two groups.

Let $\mathcal{X}=\{ X_1,\ldots, X_{n_1}\}$ be a random sample from a
one-dimensional population with mean $\mu_1$ and variance $\sigma_{1
}^2$, and let $\mathcal{Y}=\{Y_1,\ldots, Y_{n_2}\}$ be a random sample
from another one-dimensional population with mean $\mu_2$ and variance
$\sigma_2^2$ independent of $\mathcal{X}$. The two-sample $t$-statistic
is defined as
\[
\wh{T}_{\bar{n}} = \frac{\bar{X}-\bar{Y}}{\sqrt{ \wh{\sigma}_{1}^2
n_1^{-1} + \wh
{\sigma}_{2}^2 n_2^{-1} }},\label{student.t}
\]
where $\bar{n}=(n_1, n_2)$, $\bar{X}={n}^{-1}_1\sum_{i=1}^{n_1}X_i$,
$\bar{Y}={n}^{-1}_2 \sum_{j=1}^{n_2}Y_j$ and
\[
\wh{\sigma}_{1}^2 = \frac{1}{n_1-1} \sum
_{i=1}^{n_1}(X_i -\bar {X})^2,\qquad
\wh{\sigma }_{2}^2 = \frac{1}{n_2-1} \sum
_{j=1}^{n_2}(Y_j - \bar{Y})^2.
\]
The following result is a direct consequence of Theorem~\ref{thmm.1}.

\begin{theorem} \label{cor.1}
Assume that $\mu_1=\mu_2$, and $\e(|X_1|^p) < \infty, \e(|Y_1|^p) <
\infty$ for some $2< p \leq3$. Then there exist absolute constants $C,
c>0$ such that
\[
\frac{\p( \wh{T}_{\bar{n}} \geq x)}{1-\Phi(x)} = 1+O(1 ) (1 + x)^p \sum
_{\ell
=1}^2 ( v_{\ell, p}/ \sigma_{\ell}
)^p n_{\ell}^{1-p/2 }
\]
holds uniformly for $0 \leq x \leq c \min_{\ell=1, 2} \{ (\sigma
_{\ell} / v_{\ell,p}) n_{\ell}^{1/2 - 1/p} \}$, where $|O(1)|\leq C$
and $v_{1,p}= \{ \e(|X_1 - \mu_1 |^p) \}^{1/p}$, $v_{2,p}= \{ \e(|Y_1
-\mu_2|^p )\}^{1/p}$.
\end{theorem}

Motivated by a series of recent studies on the effectiveness and
accuracy of multiple-hypothesis testing using $t$-tests, we investigate
whether a higher order expansion of the relative error, as in
Theorem~1.2 of \citet{Wang2005} for self-normalized sums, holds for the
two-sample $t$-statistic,
so that one can use bootstrap calibration to correct skewness [\citet
{FanHallYao2007,DelaigleHallJin2011}] or study power properties against
sparse alternatives [\citet{WangHall2009}]. The following theorem gives a
refined Cram\'er-type moderate deviation result for $\wh{T}_{\bar{n}} $,
whose proof is placed in the supplementary material [\citet{supp}].

\begin{theorem} \label{thmm.t-stat}
Assume that $\mu_1=\mu_2$. Let $\gamma_1 = \e\{( X_1 -\mu_1)^3\}$ and
$\gamma_2 = \e\{(Y_1 -\mu_2)^3\}$ be the third central moment of
$X_1$ and $Y_1$, respectively. Moreover, assume that
$\e(|X_1|^p) < \infty, \e(|Y_1|^p)<\infty$ for some $3< p \leq4$. Then
%
\begin{eqnarray}\label{cmd.t}
\frac{\p( \wh{T}_{\bar{n}} \geq x)}{1-\Phi(x)} & = &\exp \biggl\{ - \frac{ \gamma_1 n_1^{-2} - \gamma_2 n_2^{-2} }{ 3
(\sigma_{1}^2 n_1^{-1} + \sigma_{2}^2 n_2^{-1} )^{3/2}} x^3 \biggr
\}
\nonumber
\\[-8pt]
\\[-8pt]
\nonumber
& &{}  \times \Biggl[ 1+O(1 ) \sum_{\ell=1}^2
\biggl\{ \frac{
v_{\ell,3}^3 (1 + x)}{ \sigma_\ell^3 {n}^{1/2}_\ell} + \frac
{v_{\ell
,p}^p (1 + x)^p }{ \sigma_{\ell}^p   n_{\ell}^{ p/2-1}} \biggr\} \Biggr]
\nonumber
\end{eqnarray}
holds uniformly for
%
\begin{equation}
0 \leq x \leq c \min_{\ell=1, 2} \min \bigl\{ (
\sigma_\ell/v_{\ell,
3})^3 n_\ell^{1/2}
, (\sigma_\ell/v_{\ell, p}) n_{\ell}^{1/2-1/p}
\bigr\}, \label{x.range}
\end{equation}
where $|O(1)|\leq C$ and for every $q \geq1$, $v_{1, q}= \{ \e
(|X_1 - \mu_1 |^q )\}^{1/q}$, $v_{2, q}= \{\e(|Y_1 -\mu_2|^q) \}^{1/q}$.
\end{theorem}

A refined Cram\'er-type moderate deviation theorem for the one-sample
$t$-statistic was established in \citet{Wang2011}, which to our
knowledge, is the best result for the $t$-statistic known up to date,
or equivalently, self-normalized sums.

\subsubsection{More examples of two-sample $U$-statistics}
\label{examples}

Beyond the two-sample $t$-statistic, we enumerate three more well-known
two-sample $U$-statistics and refer to \citet{NikitinPonikarov2006} for
more examples. Let $\mathcal{X}=\{X_1,\ldots, X_{n_1}\}$ and
$\mathcal
{Y}=\{Y_1, \ldots, Y_{n_2}\}$ be two independent random samples from
population distributions $P$ and $Q$, respectively.

\begin{example}[(The Mann--Whitney test statistic)]
The kernel $h$ is of order $(s_1, s_2)=(1, 1)$, defined as
\[
h(x; y) = I\{x\leq y\}-1/2  \qquad\mbox{with } \theta= \p(X_1\leq
Y_1)-1/2,
\]
and in view of \eqref{h1h2},
\[
h_1(x)=1/2 - G(x),\qquad h_2(y)=F(y)-1/2.
\]
In particular, if $F\equiv G$, we have $\sigma_1^2=\sigma_2^2 =1/12$.
\end{example}

\begin{example}[(The Lehmann statistic)] The kernel $h$ is of order
$(s_1, s_2)=(2,2)$, defined as
\[
h(x_1, x_2; y_1, y_2 ) = I\bigl\{
|x_1- x_2| \leq|y_1 - y_2| \bigr\} -
1/2
\]
with $ \theta= \p( |X_1- X_2| \leq|Y_1 - Y_2|)-1/2$. Then under
$H_0: \theta=0$, $\e\{ h(X_1, X_2;\break Y_1,  Y_2)\}=0$, and
\[
h_1(x)= G(x)\bigl\{1-G(x)\bigr\}-1/6,\qquad h_2(y) = F(y)\bigl
\{F(y)-1\bigr\} +1/6.
\]
In particular, if $F\equiv G$, then $\sigma_1^2 = \sigma_2^2 = 1/180$.
\end{example}

\begin{example}[(The Kochar statistic)] The Kochar statistic was
constructed by \citet{Kochar1979} to test if the two hazard failure
rates are different. Denote by $\mathcal{F}$ the class of all
absolutely continuous cumulative distribution functions (CDF) $F(\cdot
)$ satisfying $F(0)=0$. For two arbitrary CDF's $F, G \in\mathcal{F}$,
and let $f=F'$, $g=G'$ be their densities. Thus, the hazard failure
rates are defined by
\[
r_F(t) = \frac{f(t)}{1-F(t)},\qquad   r_G(t) =
\frac{g(t)}{1-G(t)},
\]
as long as both $1-F(t)$ and $1-G(t)$ are positive.
\citet{Kochar1979} considered the problem of testing the null hypothesis
$H_0 : r_F(t) = r_G(t)$ against the alternative $H_1: r_F(t) \leq
r_G(t), t\geq0$ with strict inequality over a set of nonzero
measures. Observe that $H_1$ holds if and only if $\de(s,t ) =
\bar{F}(s)\bar{G}(t) - \bar{F}(t) \bar{G}(s) \geq0$ for $s \geq t
\geq
0$ with strict inequality over a set of nonzero measures, where $\bar
{F}(\cdot) := 1-F(\cdot)$ for
any $F \in\mathcal{F}$.

Recall that $X_1,\ldots, X_{n_1}$ and $Y_1,\ldots, Y_{n_2}$
are two independent samples drawn respectively from $F$
and $G$. Following \citet{NikitinPonikarov2006}, we see that
\begin{eqnarray*}
\eta(F;G)& =& \e\bigl\{ \de(X \vee Y, X \wedge Y) \bigr\}
\\
&=& \p(Y_1 \leq Y_2 \leq X_1 \leq
X_2 ) + \p(X_1 \leq Y_2 \leq Y_2
\leq X_2 )
\\
&&{} -\p(X_1 \leq X_2 \leq Y_1 \leq
Y_2 ) -\p(Y_1 \leq X_1 \leq X_2
\leq Y_2 ).
\end{eqnarray*}
Under $H_0$, $\eta(F;G)=0$ while under $H_1$, $\eta(F;G)>0$. The
$U$-statistic with the kernel of order $(s_1,
s_2)=(2,2)$ is given by
\[
h(x_1, x_2; y_1, y_2) = I\{yyxx
\mbox{ or } xyyx \} - I\{ xxyy \mbox{ or } yxxy \}.
\]
Here, the term ``$yyxx$'' refers to $y_1\leq y_2 \leq x_1 \leq
x_2 $ and similar treatments apply to
$xyyx$, $xxyy$ and $yxxy$. Under $H_0 : r_F(t) = r_G(t)$, we have
\begin{eqnarray*}
h_1(x)&=& -4 G^3(x) /3 + 4 G^2(x) -2 G(x),\\
h_2(y) &=& 4 F^3(y)/3 -4 F^2(y) +2 F(y).
\end{eqnarray*}
In particular, if $F\equiv G$, then $\sigma_1^2 = \sigma_2^2 =8/105$.
\end{example}

\section{Multiple testing via Studentized two-sample tests}
\label{HT.sec}

Multiple-hypothesis testing occurs in a wide range of applications
including DNA microarray experiments, functional magnetic resonance
imaging analysis (fMRI) and astronomical surveys. We refer to \citet
{Dudoitvan2008} for a systematic study of the existing multiple testing
procedures. In this section, we consider multiple-hypothesis testing
based on Studentized two-sample tests and show how the theoretical
results in the previous section can be applied to these problems.

\subsection{Two-sample $t$-test}

A typical application of multiple-hypothesis testing in high dimensions
is the analysis of gene expression microarray data. To see whether each
gene in isolation behaves differently in a control group versus an
experimental group, we can apply the two-sample $t$-test. Assume that
the statistical model is given by
%
\begin{equation}
\cases{ X_{i,k} = \mu_{1k} +
\varepsilon_{i,k}, & \quad $i=1, \ldots, n_1$, \vspace*{2pt}
\cr
Y_{j,k} = \mu_{2k} + \omega_{j,k}, &\quad  $j=1,\ldots,
n_2$, } %
\label{two.mean.mod}
\end{equation}
for $k=1, \ldots, m$, where index $k$ denotes the $k$th gene, $i$ and
$j$ indicate the $i$th and $j$th array, and the constants $\mu_{1k}$
and $\mu_{2k}$, respectively, represent the mean effects for the $k$th
gene from the first and the second groups. For each $k$, $\varepsilon
_{1,k}, \ldots, \varepsilon_{n_1,k}$ (resp., $\omega_{1,k}, \ldots,
\omega_{n_2 , k}$) are independent random variables with mean zero and
variance $\sigma_{1k}^2>0$ (resp., $\sigma_{2k}^2>0$). For the $k$th
marginal test, when the population variances $\sigma_{1k}^2$ and
$\sigma
_{2k}^2$ are unequal, the two-sample $t$-statistic is most commonly
used to carry out hypothesis testing for the null $H^k_{0}:\mu_{1k}
=\mu
_{2k}$ against the alternative $H^k_{1}:\mu_{1k} \neq\mu_{2k}$.

Since the seminal work of \citet{BenjaminiHochberg1995}, the Benjamini
and Hochberg (B--H) procedure has become a popular technique in
microarray data analysis for gene selection, which along with many
other procedures depend on $p$-values that often need to be estimated.
To control certain simultaneous errors, it has been shown that using
approximated $p$-values is asymptotically equivalent to using the true
$p$-values for controlling the $k$-familywise error rate ($k$-FWER) and
false discovery rate (FDR). See, for example, \citet{KosorokMa2007},
\citet{FanHallYao2007} and \citet{LiuShao2010} for one-sample tests.
\citet
{CaoKosorok2011} proposed an alternative method to control $k$-FWER and
FDR in both large-scale one- and two-sample $t$-tests. A common thread
among the aforementioned literature is that theoretically for the
methods to work in controlling FDR at a given level, the number of
features $m$ and the sample size $n$ should satisfy $\log m = o(n^{1/3})$.

Recently, \citet{LiuShao2014} proposed a regularized bootstrap
correction method for multiple one-sample $t$-tests so that the
constraint on $m$ may be relaxed to $\log m= o(n^{1/2})$ under less
stringent moment conditions as assumed in \citet{FanHallYao2007} and
\citet{DelaigleHallJin2011}. Using Theorem~\ref{thmm.t-stat}, we show
that the constraint on $m$ in large scale two-sample $t$-tests can be
relaxed to $\log m = o(n^{1/2})$ as well. This provides theoretical
justification of the effectiveness of the bootstrap method which is
frequently used for skewness correction.

To illustrate the main idea, here we restrict our attention to the
special case in which the observations are independent. Indeed, when
test statistics are correlated, false discovery control becomes very
challenging under arbitrary dependence. Various dependence structures
have been considered in the literature. See, for example, \citet
{BenjaminiYekutieli2001}, \citet{StoreyTaylorSiegmund2004}, \citet
{FerreiraZwinderman2006}, \citet{LeekStorey2008}, \citet
{FriguetKloaregCauseur2009} and \citet{FanHanGu2012}, among others. For
completeness, we generalize the results to the dependent case in
Section~\ref{se:depen}.

\subsubsection{Normal calibration and phase transition}\label{se:norm}

Consider the large-scale significance testing problem:
\[
H^k_0 : \mu_{1k}= \mu_{2k}
\quad\mbox{versus}\quad H^k_1 : \mu_{1k}\neq
\mu_{2k}, \qquad 1\leq k \leq m.
\]
Let $V$ and $R$ denote, respectively, the number of false rejections
and the number of total rejections. The well-known false discovery
proportion (FDP) is defined as the ratio $\mathrm{FDP}= V/\max(1, R)$, and FDR
is the expected FDP, that is, $\e\{V/\max(1,R)\}$. \citet
{BenjaminiHochberg1995} proposed a\break distribution-free method for
choosing a $p$-value threshold that controls the FDR at a prespecified
level where $0<\alpha<1$. For $k=1, \ldots, m$, let $p_k$ be the
marginal $p$-value of the $k$th test, and let $p_{(1)}\leq\cdots\leq
p_{(m)}$ be the order statistics of $p_1, \ldots, p_m$. For a
predetermined control level $\alpha\in(0,1)$, the B--H procedure
rejects hypotheses for which $p_k \leq p_{( \hat{k})}$, where
%
\begin{equation}
\hat{k} = \max \biggl\{ 0\leq k\leq m : p_{(k)} \leq\frac{\alpha k}{m}
\biggr\} \label{p-threshold}
\end{equation}
with $p_{(0)}=0$.

In microarray analysis, two-sample $t$-tests are often used to identify
differentially expressed genes between two groups. Let
\[
T_k = \frac{\bar{X}_k - \bar{Y}_k}{\sqrt{\wh{\sigma}_{1k}^2 n_1^{-1}+
\wh
{\sigma}_{2k}^2 n_2^{-1}}},\qquad k=1, \ldots, m,
\]
where $ \bar{X}_k = n_1^{-1}\sum_{i=1}^{n_1}X_{i,k}$, $\bar{Y}_k =
n_2^{-1}\sum_{j=1}^{n_2}
Y_{j,k}$ and
\[
\wh{\sigma}_{1k}^2=\frac{1}{n_1-1} \sum
_{i=1}^{n_1}(X_{i,k}-\bar
{X}_k)^2, \qquad  \wh {\sigma}_{2k}^2=
\frac{1}{n_2-1} \sum_{j=1}^{n_2}(Y_{j,k}-
\bar{Y}_k)^2.
\]
Here and below, $\{X_{i,1}, \ldots, X_{i,m}\}_{i=1}^{n_1}$ and $\{
Y_{j,1}, \ldots, Y_{j,m}\}_{j=1}^{n_2}$ are independent random samples
from $\{X_1, \ldots, X_m\}$ and $\{Y_1, \ldots, Y_m\}$, respectively,
generated according to model (\ref{two.mean.mod}), which are usually
non-Gaussian in practice. Moreover, assume that the sample sizes of the
two samples are of the same order, that is, $n_1 \asymp n_2$.

Before stating the main results, we first introduce a number of
notation. Set $\mathcal{H}_0=\{ 1\leq k\leq m : \mu_{1k}= \mu_{2k}\}$,
let $m_0 = \# \mathcal{H}_0 $ denote the number of true null hypotheses
and $m_1=m-m_0$. Both $m=m(n_1, n_2)$ and $m_0=m_0(n_1,n_2)$ are
allowed to grow as $n=n_1 \wedge n_2$ increases. We assume that
\[
\lim_{n \rightarrow\infty} \frac{m_0}{m} = \pi_0 \in(0,1].
\]
In line with the notation used in Section~\ref{cmd.U}, set
\begin{eqnarray*}
\sigma_{1k}^2 &= &\var(X_k),\qquad
\sigma_{2k}^2=\var(Y_k) ,\\
 \gamma
_{1k}& =& \e\bigl\{(X_k-\mu_{1k})^3
\bigr\} ,\qquad \gamma_{2k} = \e\bigl\{(Y_k-\mu
_{2k})^3\bigr\}
\end{eqnarray*}
and $\sigma_{\bar{n},k}^2 = \sigma_{1k}^2 n_1^{-1} + \sigma_{2k}^2
n_2^{-1}$. Throughout this subsection, we focus on the normal
calibration and let $\wh{p}_k=2-2 \Phi(|T_k|)$, where $\Phi(\cdot)$ is
the standard normal distribution function. Indeed, the exact null
distribution of $T_k$ and thus the true $p$-values are unknown without
the normality assumption.

\begin{theorem} \label{app.thmm.1}
Assume that $\{X_1, \ldots, X_m, Y_1, \ldots, Y_m\}$ are independent
nondegenerate random variables; $n_1 \asymp n_2$, $m=m(n_1,
n_2)\rightarrow\infty$ and $\log m = o(n^{1/2})$ as $n=n_1 \wedge n_2
\rightarrow\infty$. For independent random samples $\{X_{i,1}, \ldots,\break
X_{i,m}\}_{i=1}^{n_1}$ and $\{Y_{j,1}, \ldots, Y_{j,m}\}_{j=1}^{n_2}$,
suppose that
%
\begin{equation}
\min_{1\leq k\leq m}\min( \sigma_{1k}, \sigma_{2k}
) \geq c>0,\qquad \max_{1\leq k\leq m}\max \bigl\{ \e\bigl(
\xi_k^4\bigr) , \e\bigl(\eta _k^4
\bigr) \bigr\} \leq C < \infty\label{4th.mc}
\end{equation}
for some constants $C$
and $c$, where $\xi_k = \sigma_{1k}^{-1} (X_k-\mu_{1k}) $ and $\eta
_k =
\sigma_{2k}^{-1} (Y_k-\mu_{2k})$. Moreover, assume that
%
\begin{equation}
\# \bigl\{ 1\leq k\leq m : |\mu_{1k}-\mu_{2k}| \geq4 (\log
m)^{1/2} \sigma_{\bar{n}, k} \bigr\} \rightarrow\infty\label{sign}
\end{equation}
as $n\rightarrow\infty$, and let
%
\begin{equation}
c_0 = \liminf_{n, m\rightarrow\infty} \biggl\{ \frac
{{n}^{1/2}}{m_0}
\sum_{k \in\mathcal{H}_{0}} \sigma_{\bar{n}, k}^{-3} \bigl|
\gamma_{1k} n_1^{-2} - \gamma_{2k}
n_2^{-2} \bigr| \biggr\}. \label{skew.parameter}
\end{equation}
\begin{longlist}[(iii)]
\item[(i)] Suppose that $\log m= o(n^{1/3})$. Then as $n \to
\infty
$, $\mathrm{ FDP}_\Phi\rightarrow^{ P} \alpha\pi_0$ and $\mathrm{ FDR}_\Phi
\to\alpha\pi_0$.

\item[(ii)] Suppose that $c_0>0$, $\log m \geq c_1 n^{1/3}$ for
some $c_1>0$ and that $\log m_1 =o(n^{1/3})$. Then there exists some
constant $\beta\in(\alpha, 1]$ such that
\[
\lim_{n \to\infty} \p( \mathrm{ FDP}_\Phi\geq\beta) = 1\quad \mbox{and}\quad \liminf_{n \to\infty} \mathrm{ FDR}_\Phi\geq\beta.
\]

\item[(iii)] Suppose that $c_0>0$, $( \log m ) /n^{1/3}
\rightarrow\infty$ and $ \log m_1 = o( n^{1/3} )$. Then as $n \to
\infty$, $\mathrm{ FDP}_{\Phi} \rightarrow^{ P} 1$ and $\mathrm{ FDR}_\Phi
\to1$.
\end{longlist}
Here, $\mathrm{ FDR}_{\Phi}$ and $\mathrm{ FDP}_\Phi$ denote, respectively, the
$\mathrm{ FDR}$ and the $\mathrm{ FDP}$ of the B--H procedure with $p_k$ replaced
by $\wh{p}_k$ in (\ref{p-threshold}).
\end{theorem}

Together, conclusions (i) and (ii) of Theorem~\ref{app.thmm.1} indicate
that the number of simultaneous tests can be as large as $\exp\{
o(n^{1/3})\}$ before the normal calibration becomes inaccurate. In
particular, when $n_1=n_2=n$, the skewness parameter $c_0$ given in
(\ref{skew.parameter}) reduces to
\[
c_0 = \liminf_{ m \rightarrow\infty} \biggl\{ \frac{1}{m_0}
\sum_{ k
\in
\mathcal{H}_{0}} \frac{|\gamma_{1k} - \gamma_{2k} | }{(\sigma
_{1k}^2 +
\sigma_{2k}^2)^{ 3/2} } \biggr\}.
\]
As noted in \citet{LiuShao2014}, the limiting behavior of the
$\mathrm{FDR}_{\Phi
}$ varies in different regimes and exhibits interesting phase
transition phenomena as the dimension $m$ grows as a function of $(n_1,
n_2)$. The average of skewness $c_0$ plays a crucial role. It is also
worth noting that conclusions (ii) and (iii) hold under the scenario
$\pi_0=1$, that is, $m_1=o(m)$. This corresponds to the sparse settings
in applications such as gene detections. Under finite 4th moments of
$X_k$ and $Y_k$, the robustness of two-sample $t$-tests and the
accuracy of normal calibration in the FDR/FDP control have been
investigated in \citet{CaoKosorok2011} when $m_1/m \to\pi_1 \in(0,1)$.
This corresponds to the relatively dense setting, and the sparse case
that we considered above is not covered.

\subsubsection{Bootstrap calibration and regularized bootstrap correction}
\label{reg.bootstrap}

In this subsection, we first use the conventional bootstrap calibration
to improve the accuracy of FDR control based on the fact that the
bootstrap approximation removes the skewness term that determines
first-order inaccuracies of the standard normal approximation. However,
the validity of bootstrap approximation requires the underlying
distribution to be very light tailed, which does not seem realistic in
real data applications. As pointed in the literature of gene study,
many gene data are commonly recognized to have heavy tails which
violates the assumption on underlying distribution used to make
conventional bootstrap approximation work. Recently, \citet{LiuShao2014}
proposed a regularized bootstrap method that is shown to be more robust
against the heavy tailedness of the underlying distribution and the
dimension $m$ is allowed to be as large as $\exp\{ o(n^{1/2})\}$.

Let $\mathcal{X}_{k,b}^{\dagger}=\{X^{\dagger}_{1,k,b}, \ldots,
X^{\dagger}_{n_1,k, b}\}$, $\mathcal{Y}_{k,b}^{\dagger}=\{Y^{\dagger
}_{1,k,b}, \ldots, Y^{\dagger}_{n_2,k,b}\}$, $b=1, \ldots, B$, denote
bootstrap samples drawn independently and uniformly, with replacement,
from $\mathcal{X}_k=\{X_{1,k}, \ldots, X_{n_1,k}\}$ and $\mathcal
{Y}_k=\{Y_{1,k}, \ldots, Y_{n_2,k}\}$, respectively. Let $T^{\dagger
}_{k,b}$ be the two-sample $t$-statistic constructed from $\{
X_{1,k,b}^{\dagger}- \bar{X}_k , \ldots,\break X^{\dagger}_{n_1, k,b} -
\bar
{X}_k \}$ and $\{ Y^{\dagger}_{1,k,b}-\bar{Y}_k , \ldots, Y^{\dagger
}_{n_2,k,b} -\bar{Y}_k \}$. Following \citet{LiuShao2014}, we use the
following empirical distribution:
\[
F^{\dagger}_{m,B}(t) = \frac{1}{m B } \sum
_{k=1}^m \sum_{b=1}^B
I\bigl\{ \bigl|T^{\dagger}_{k,b}\bigr| \geq t\bigr\}
\]
to approximate the null distribution, and thus the estimated $p$-values
are given by $\wh{p}_{k, \mathrm{ B} }=F^{\dagger}_{m,B}(|T_k|)$.
Respectively, $\mathrm{FDP}_{\mathrm{ B}}$ and $\mathrm{FDR}_{\mathrm{ B}}$ denote the FDP and
the FDR of the B--H procedure with $p_k$ replaced by $\wh{p}_{k, \mathrm{
B}}$ in (\ref{p-threshold}).

The following result shows that the bootstrap calibration is accurate
provided $\log m$ increases at a strictly slower rate than $(n_1 \wedge
n_2)^{1/2}$, and the underlying distribution has sub-Gaussian tails.

\begin{theorem} \label{app.thmm.2}
Assume the conditions in Theorem~\ref{app.thmm.1} hold and that
\[
\max_{1\leq k\leq m} \max \bigl\{ \e\bigl(e^{t_0 \xi_k^2}\bigr) , \e
\bigl( e^{t_0
\eta
_k^2}\bigr) \bigr\} \leq C < \infty
\]
for some constants $t_0, C>0$.
\begin{longlist}[(ii)]
\item[(i)] Suppose that $\log m = o(n^{1/3})$. Then as $n \to
\infty$, $\mathrm{ FDP}_\mathrm{ B} \rightarrow^{ P} \alpha\pi_0$ and $\mathrm{
FDR}_\mathrm{ B} \to\alpha\pi_0$.

\item[(ii)] Suppose that $\log m = o(n^{1/2})$ and $m_1 \leq
m^{\rho}$ for some $\rho\in(0,1)$. Then as $n \to\infty$, $\mathrm{
FDP}_\mathrm{ B} \rightarrow^{ P} \alpha$ and $\mathrm{ FDR}_\mathrm{ B} \to
\alpha$.
\end{longlist}
\end{theorem}

The sub-Gaussian condition in Theorem~\ref{app.thmm.2} is quite
stringent in practice, whereas it can hardly be weakened in general
when the bootstrap method is applied. In the context of family-wise
error rate control, \citet{FanHallYao2007} proved that the bootstrap
calibration is accurate if the observed data are bounded and $\log m =
o(n^{1/2})$. The regularized bootstrap method, however, adopts the very
similar idea of the trimmed estimators and is a two-step procedure that
combines the truncation technique and the bootstrap method.

First, define the trimmed samples
\[
\widehat{X}_{i,k} = X_{i,k}I\bigl\{|X_{i,k}| \leq
\la_{1k}\bigr\},  \qquad \widehat{Y}_{j,k} = Y_{i,k}I
\bigl\{|Y_{j,k}| \leq\la_{2k}\bigr\}
\]
for $i=1, \ldots, n_1$, $j=1, \ldots, n_2$, where $\la_{1k}$ and
$\la_{2k}$ are regularized parameters to be specified. Let $\wh
{\mathcal
{X}}_{k,b}^{\dagger}=\{\wh{X}^{\dagger}_{1, k,b }, \ldots, \wh
{X}^{\dagger}_{n_1, k,b}\}$ and $\wh{\mathcal{Y}}_{k,b}^{\dagger
}=\{ \wh
{Y}^{\dagger}_{1, k,b }, \ldots, \wh{Y}^{\dagger}_{n_2, k,b }\}$, $b=1,
\ldots, B$, be the corresponding bootstrap samples drawn by sampling
randomly, with replacement, from
\[
\wh{\mathcal{X}}_k =\{\wh{X}_{1,k}, \ldots,
\wh{X}_{ n_1, k}\}\quad \mbox{and}\quad \wh{\mathcal{Y}}_k =\{
\wh{Y}_{1,k}, \ldots, \wh{Y}_{
n_2 , k}\},
\]
respectively. Next, let $\wh{T}^{\dagger}_{k,b}$ be the two-sample
$t$-test statistic constructed from $\{ \wh{X}_{1, k,b }^{\dagger}-
n_1^{-1}\sum_{i=1}^{n_1}\wh{X}_{i,k}, \ldots, \wh{X}^{\dagger
}_{n_1, k, b} -
n_1^{-1}\sum_{i=1}^{n_1}\wh{X}_{i,k} \}$ and $\{ \wh{Y}^{\dagger
}_{1, k,b }-
n_2^{-1} \times \break \sum_{j=1}^{n_2}\wh{Y}_{j,k} , \ldots, \wh{Y}^{\dagger
}_{n_2, k, b}
-n_2^{-1}\sum_{j=1}^{n_2}\wh{Y}_{j,k} \}$. As in the previous
procedure, define the
estimated
$p$-values by
\[
\wh{p}_{k,\mathrm{ RB}} = \wh{F}^{\dagger}_{m, \mathrm{ RB}}\bigl(|T_k|\bigr)
\qquad\mbox{with } \wh{F}^{\dagger}_{m, \mathrm{ RB} }(t)= \frac{1}{ m B } \sum
_{k=1}^m \sum
_{b=1}^B I\bigl\{\bigl|\wh{T}^{\dagger}_{k,b}\bigr|
\geq t\bigr\}.
\]
Let $\mathrm{FDP}_{\mathrm{ RB}}$ and $\mathrm{FDR}_{\mathrm{ RB}}$ denote the FDP and the FDR,
respectively, of the B--H procedure with $p_k$ replaced by $\wh{p}_{ k,
\mathrm{ RB}}$ in (\ref{p-threshold}).

\begin{theorem} \label{app.thmm.3}
Assume the conditions in Theorem~\ref{app.thmm.1} hold and that
%
\begin{equation}
\max_{1\leq k \leq m} \max \bigl\{ \e\bigl(|X_k|^6
\bigr) , \e\bigl(|Y_k|^6\bigr) \bigr\} \leq C < \infty.
\label{6th.moment}
\end{equation}
The regularized parameters $(\lambda_{1k}, \lambda_{2k})$ are such that
%
\begin{equation}
\la_{ 1k } \asymp \biggl( \frac{n_1}{\log m} \biggr)^{1/6}
  \quad\mbox {and} \quad \la_{ 2k } \asymp \biggl( \frac{n_2}{\log m}
\biggr)^{1/6}. \label{reg.order}
\end{equation}
\begin{longlist}[(ii)]
\item[(i)] Suppose that $\log m= o(n^{1/3})$. Then as $n \to
\infty
$, $\mathrm{ FDP}_{\mathrm{ RB}} \rightarrow^{ P} \alpha\pi_0$ and $\mathrm{
FDR}_{\mathrm{ RB}} \to\alpha\pi_0$.

\item[(ii)] Suppose that $\log m= o(n^{1/2})$ and $m_1 \leq
m^{\rho
}$ for some $\rho\in(0,1)$. Then as $n \to\infty$, $\mathrm{ FDP}_\mathrm{
RB} \rightarrow^{ P} \alpha$ and $\mathrm{ FDR}_\mathrm{ RB} \to\alpha$.
\end{longlist}
\end{theorem}

In view of Theorem~\ref{app.thmm.3}, the regularized bootstrap
approximation is valid under mild moment conditions that are
significantly weaker than those required for the bootstrap method to
work theoretically. The numerical performance will be investigated in
Section~\ref{sim.sec}. To highlight the main idea, a self-contained
proof of Theorem~\ref{app.thmm.1} is given in the supplementary material [\citet{supp}]. The proofs of Theorems \ref{app.thmm.2} and \ref{app.thmm.3}
are based on straightforward extensions of Theorems 2.2 and 3.1 in
\citet
{LiuShao2014}, and thus are omitted.

\subsubsection{FDR control under dependence}\label{se:depen}

In this section, we generalize the results in previous sections to the
dependence case. Write $\varrho=n_1/n_2$. For every $k, \ell=1,\ldots,
m$, let $\sigma_k^2 = \sigma_{1k}^2 + \varrho\sigma_{2k}^2$ and define
%
\begin{equation}
r_{k \ell} = (\sigma_k \sigma_{\ell}
)^{-1} \bigl\{ \cov(X_{k}, X_{\ell
}) + \varrho
\cov(Y_k, Y_{\ell}) \bigr\}, \label{corr.def}
\end{equation}
which characterizes the dependence between $(X_k, Y_k)$ and
$(X_\ell, Y_\ell)$. Particularly, when $n_1=n_2$ and $\sigma
_{1k}^2=\sigma_{2k}^2$, we see that $r_{k \ell}= \frac{1}{2} \{
\corr
(X_k, X_{\ell}) + \corr(Y_k, Y_{\ell}) \}$. In this subsection, we
impose the following conditions on the dependence structure of $\bX=
(X_1,\ldots, X_m)^{\mathrm{T}}$ and $\bY= (Y_1, \ldots,
Y_m)^{\mathrm{T}}$.
\begin{longlist}[(D1)]
\item[(D1)] There exist constants $0<r<1$, $0<\rho< (1-r)/(1+r) $ and
$b_1 >0$ such that
\[
\max_{1\leq k\neq\ell\leq m} | r_{k \ell} | \leq r \quad\mbox{and}\quad \max
_{1\leq k \leq m } s_k(m) \leq b_1
m^\rho,
\]
where for $k =1,\ldots, m$,
\begin{eqnarray*}
s_k(m)& =& \bigl\{ 1\leq\ell\leq m : \corr(X_k,
X_{\ell}) \geq (\log m)^{-2-\gamma}
\nonumber
\\
& &      \mbox{or }    \corr (Y_k,
Y_{\ell}) \geq(\log m)^{-2-\gamma} \bigr\}
\end{eqnarray*}
for some $\gamma>0$.

\item[(D2)] There exist constants $0<r<1$, $0<\rho< (1-r)/(1+r) $ and
$b_1 >0$ such that $\max_{1\leq k\neq\ell\leq m} | r_{k \ell} |
\leq r$ and for each $X_k$, the number of variables $X_\ell$ that are
dependent of $X_k$ is less than $b_1 m^\rho$.
\end{longlist}

The assumption $\max_{1\leq k\neq\ell\leq m} | r_{k \ell} | \leq r$
for some $0<r<1$ imposes a constraint on the magnitudes of the
correlations, which is natural in the sense that the correlation matrix
$\mathbf{R}=(r_{k \ell})_{1\leq k, \ell\leq m }$ is singular if
$\max_{1\leq k\neq\ell\leq m} | r_{k \ell} | = 1$. Under condition (D1),
each $(X_k,Y_k)$ is allowed to be ``moderately'' correlated with at
most as many as $O(m^\rho)$ other vectors. Condition (D2) enforces a
local dependence structure on the data, saying that each vector is
dependent with at most as many as $O(m^\rho)$ other random vectors and
independent of the remaining ones. The following theorem extends the
results in previous sections to the dependence case. Its proof is
placed in the supplementary material [\citet{supp}].

\begin{theorem} \label{FDR.dependence}
Assume that either condition \textup{(D1)} holds with $\log m = O(n^{1/8})$
or condition \textup{(D2)} holds with $\log m = o(n^{1/3})$.
\begin{longlist}[(ii)]
\item[(i)] Suppose that (\ref{4th.mc}) and (\ref{sign})
are satisfied. Then as $n\to\infty$, $\mathrm{ FDP}_\Phi\rightarrow^{ P}
\alpha\pi_0$ and $\mathrm{ FDR}_\Phi\to\alpha\pi_0$.

\item[(ii)] Suppose that (\ref{4th.mc}), (\ref
{6th.moment}) and (\ref{reg.order}) are satisfied. Then as
$n\to
\infty$,  $\mathrm{ FDP}_\mathrm{ RB} \rightarrow^{ P} \alpha\pi_0$ and
$\mathrm{
FDR}_\mathrm{ RB} \to\alpha\pi_0$.
\end{longlist}
In particular, assume that condition \textup{(D2)} holds with $\log m =
o(n^{1/2})$ and $m_1 \leq m^c$ for some $0<c<1$. Then as $n\to\infty$,
$\mathrm{ FDP}_\mathrm{ RB} \rightarrow^{ P} \alpha\pi_0$ and $\mathrm{
FDR}_\mathrm{
RB} \to\alpha\pi_0$.
\end{theorem}


\subsection{Studentized Mann--Whitney test}
\label{app.MW}

Let $\mathcal{X}=\{X_1, \ldots, X_{n_1}\}$ and $\mathcal{Y}=\{Y_1,
\ldots, Y_{n_2}\}$ be two independent random samples from distributions
$F$ and $G$, respectively. Let $\theta=\p(X\leq Y)-1/2$. Consider the
null hypothesis $H_0: \theta=0$ against the one-sided alternative $H_1:
\theta>0$. This problem arises in many applications including testing
whether the physiological performance of an active drug is better than
that under the control treatment, and testing the effects of a policy,
such as unemployment insurance or a vocational training program, on the
level of unemployment.

The Mann--Whitney (M--W) test [\citet{MannWhitney1947}], also known as the
two-sample Wilcoxon test [\citet{Wilcoxon1945}], is prevalently used for
testing equality of means or medians, and serves as a nonparametric
alternative to the two-sample $t$-test. The corresponding test
statistic is given by
%
\begin{equation}
U_{\bar{n}} = \frac{1}{n_1 n_2} \sum_{i=1}^{n_1}
\sum_{j=1}^{n_2}I\{ X_i \leq
Y_j \}, \qquad \bar{n}=(n_1, n_2).
\label{M-T.stat}
\end{equation}

The M--W test is widely used in a wide range of fields including
statistics, economics and biomedicine, due to its good efficiency and
robustness against parametric assumptions. Over one-third of the
articles published in \textit{Experimental Economics} use the
Mann--Whitney test and \citet{Okeh2009} reported that thirty percent of
the articles in five biomedical journals published in 2004 used the
Mann--Whitney test. For example, using the M--W $U$ test, \citet
{CharnessGneezy2009} developed an experiment to test the conjecture
that financial incentives help to foster good habits. They recorded
seven biometric measures (weight, body fat percentage, waist size,
etc.) of each participant before and after the experiment to assess the
improvements across treatments. Although the M--W test was originally
introduced as a rank statistic to test if the distributions of two
related samples are identical, it has been prevalently used for testing
equality of medians or means, sometimes as an alternative to the
two-sample $t$-test.

It was argued and formally examined recently in \citet{ChungRomano2011}
that the M--W test has generally been misused across disciplines. In
fact, the M--W test is only valid if the underlying distributions of the
two groups are identical. Nevertheless, when the purpose is to test the
equality of distributions, it is recommended to use a statistic, such
as the Kolmogorov--Smirnov or the Cram\'er--von Mises statistic, that
captures the discrepancies of the entire distributions rather than an
individual parameter. More specifically, because the M--W test only
recognizes deviation from $\theta=0$, it does not have much power in
detecting overall distributional discrepancies. Alternatively, the M--W
test is frequently used to test the equality of medians. However, \citet
{ChungRomano2013} presented evidence that this is another improper
application of the M--W test and suggested to use the Studentized median test.

Even when the M--W test is appropriately applied for testing $H_0:
\theta
=0$, the asymptotic variance depends on the underlying distributions,
unless the two population distributions are identical. As \citet
{HallWilson1991} pointed out, the application of resampling to pivotal
statistics has better asymptotic properties in the sense that the rate
of convergence of the actual significance level to the nominal
significance level is more rapid when the pivotal statistics are
resampled. Therefore, it is natural to use the Studentized
Mann--Whitney test, which is asymptotic pivotal.

Let
%
\begin{equation}
\wh{U}_{\bar{n}} = \wh{\sigma}_{\bar{n}}^{-1}(U_{\bar{n}}-
1/2) \label{stu.wilcoxon}
\end{equation}
denote the Studentized test statistic for $U_{\bar{n}}$ as in \eqref
{M-T.stat}, where $\wh{\sigma}_{\bar{n}}^2= \wh{\sigma}_1^2
n_1^{-1} +
\wh
{\sigma}_2^2 n_2^{-1}$,
\[
\wh{\sigma}_1^2 = \frac{1}{n_1-1}\sum
_{i=1}^{n_1} \Biggl( q_i -
\frac{1}{n_1}\sum_{i=1}^{n_1}q_i
\Biggr)^2 ,  \qquad  \wh{\sigma}_2^2 =
\frac{1}{n_2-1}\sum_{j=1}^{n_2} \Biggl(
p_j - \frac{1}{n_2}\sum_{j=1}^{n_2}p_j
\Biggr)^2
\]
with $q_i = n_2^{-1}\sum_{j=1}^{n_2}I\{Y_j < X_i\}$ and $p_j=
n_1^{-1}\sum_{i=1}^{n_1}I\{ X_i \leq Y_j\}$.

When dealing with samples from a large number of geographical regions
(suburbs, states, health service areas, etc.), one may need to make
many statistical inferences simultaneously. Suppose we observe a family
of paired groups, that is, for $k=1,\ldots, m$, $\mathcal{X}_k = \{
X_{1,k},\ldots, X_{n_1,k}\}$, $\mathcal{Y}_k = \{Y_{1,k},\ldots,
Y_{n_2,k}\}$, where the index $k$ denotes the $k$th site. Assume that
$\mathcal{X}_k$ is drawn from $F_k$, and independently, $\mathcal{Y}_k$
is drawn from $G_k$.

For each $k=1,\ldots, m$, we test the null hypothesis $H^k_0: \theta
_k=\p(X_{1,k}\leq Y_{1,k})- 1/2 =0$ against the one-sided alternative
$H^k_1: \theta_k>0$. If $H^k_0$ is rejected, we conclude that the
treatment effect (of a drug or a policy) is acting within the $k$th
area. Define the test statistic
\[
\wh{U}_{\bar{n},k} = \wh{\sigma}_{\bar{n}, k}^{-1}(U_{\bar{n},
k}-
1/2),
\]
where $\wh{U}_{ \bar{n}, k }$ is constructed from the $k$th paired samples
according to (\ref{stu.wilcoxon}). Let
\[
F_{\bar{n},k}(t) = \p\bigl( \wh{U}_{\bar{n}, k} \leq t |
H^k_0 \bigr) \quad\mbox{and}\quad \Phi(t) = \p( Z \leq t ),
\label{p-values}
\]
where $Z$ is the standard normal random variable. Then the true
$p$-values are $p_k=1- F_{\bar{n},k}(\wh{U}_{\bar{n},k})$, and $\wh{p}_k=1-
\Phi
(\wh{U}_{\bar{n},k})$ denote the estimated $p$-values based on normal
calibration.

To identify areas where the treatment effect is acting, we can use the
B--H method to control the FDR at $\alpha$ level by rejecting the null
hypotheses indexed by $\mathcal{S}=\{ 1\leq k\leq m : \wh{p}_k \leq
\wh
{p}_{(\hat{k})} \}$, where $\hat{k} = \max\{1\leq k \leq m : \widehat{p}_{(k)}\leq \alpha k/m\}$, and $\{\wh{p}_{(k)}\}$ denote the
ordered values of $\{\wh{p}_k\}$. As before, let $\mathrm{FDR}_{\Phi}$ be the
FDR of the B--H method based on normal calibration.

Alternative to normal calibration, we can also consider bootstrap
calibration. Recall that $\mathcal{X}_{k,b}^{\dagger}=\{ X^{\dagger
}_{1,k,b}, \ldots, X^{\dagger}_{n_1,k, b} \}$ and $\mathcal
{Y}_{k,b}^{\dagger}=\{Y^{\dagger}_{1,k,b}, \ldots, Y^{\dagger
}_{n_2,k,b}\}$, $b=1, \ldots, B$, are two bootstrap samples drawn
independently and uniformly, with replacement, from $\mathcal{X}_k=\{
X_{1,k}, \ldots, X_{n_1,k}\}$ and $\mathcal{Y}_k=\{Y_{1,k}, \ldots,
Y_{n_2,k}\}$, respectively. For each $k=1,\ldots, m$, let $ \wh
{U}_{\bar{n}
,k,b}^{\dagger}$ be the bootstrapped test statistic constructed from
$\mathcal{X}_{k,b}^{\dagger}$ and $\mathcal{Y}_{k,b}^{\dagger}$, that is,
\[
\widehat{U}_{\bar{n},k,b}^\dagger=\widehat{\sigma}_{\bar
{n},k,b}^{-1}
\Biggl[U_{\bar{n},k,b}-\frac{1}{n_1n_2}\sum_{i=1}^{n_1}
\sum_{j=1}^{n_2}I\{X_{i,k}\leq
Y_{j,k}\} \Biggr],
\]
where $U_{\bar{n},k,b}$ and $\widehat{\sigma}_{\bar{n},k,b}$ are the
analogues of $U_{\bar{n}}$ given in (\ref{M-T.stat}) and $\widehat
{\sigma}_{\bar{n}}$ specified below (\ref{stu.wilcoxon}) via replacing
$X_i$ and $Y_j$ by ${X}_{i,k,b}^{\dagger}$ and ${Y}_{j,k,b}^{\dagger}$,
respectively. Using the empirical distribution function
\[
\wh G^{\dagger}_{m,B}(t) =\frac{1}{mB} \sum
_{k=1}^m \sum_{b=1}^B
I\bigl\{ \bigl| \wh{U}_{\bar{n},k,b}^{\dagger} \bigr| \leq t \bigr\},
\]
we estimate the unknown $p$-values by $\wh{p}_{k,\mathrm{ B}} = 1- \wh
G^{\dagger}_{m,B}(\wh{U}_{\bar{n},k,b}^{\dagger} ) $. For a predetermined
$\alpha\in(0,1)$, the null hypotheses indexed by $\mathcal{S}_{\mathrm{
B}} = \{ 1\leq k\leq m: \wh p_{k,\mathrm{ B}} \leq\wh p_{(\hat{k}_{\mathrm{
B}} ),\mathrm{ B}} \}$ are rejected,\vspace*{1pt} where $ \hat{k}_{\mathrm{ B}} = \max\{
0\leq k \leq m : \wh{p}_{k, \mathrm{ B}} \leq\alpha{k}/{m} \}$. Denote by
$\mathrm{FDR}_{\mathrm{ B}}$ the FDR of the B--H method based on bootstrap calibration.

Applying the general moderate deviation result (\ref{cmd}) to Studentized
Mann--Whitney statistics $\wh U_{\bar{n}, k}$ leads to the following
result. The proof is based on a straightforward adaptation of the
arguments we used in the proof of Theorem~\ref{app.thmm.1}, and hence
is omitted.

\begin{theorem} \label{MW.thmm.1}
Assume that $\{X_1, \ldots, X_m, Y_1, \ldots, Y_m\}$ are independent
random variables with continuous distribution functions $X_k \sim F_k$
and $Y_j \sim G_k$. The triplet $(n_1,n_2, m)$ is such that $n_1 \asymp
n_2$, $m=m(n_1, n_2)\rightarrow\infty$, $\log m = o(n^{1/3})$ and
$m^{-1} \# \{k=1, \ldots, m: \theta_k = 1/2\} \to\pi_0\in(0,1]$ as
$n=n_1 \wedge n_2 \rightarrow\infty$. For independent samples $\{
X_{i,1}, \ldots, X_{i,m}\}_{i=1}^{n_1}$ and $\{Y_{j,1}, \ldots,
Y_{j,m}\}_{j=1}^{n_2}$, suppose that $\min_{1\leq k\leq m}\min(
\sigma
_{1k}, \sigma_{2k} ) \geq c>0$ for some constant $c>0$ and as $n\to
\infty$,
%
\[
\# \bigl\{ 1\leq k\leq m : | \theta_k - 1/2 | \geq4 (\log
m)^{1/2} \sigma_{\bar{n}, k} \bigr\} \rightarrow\infty,
\]
where $\sigma^2_{1k}=\operatorname{ var}\{G_k(X_k)\}$, $\sigma_{2k}^2=\operatorname{
var}\{
F_k(Y_k)\}$ and $ \sigma_{\bar{n}, k}^2 =\sigma^2_{1k} n_1^{-1} +
\sigma
^2_{2k} n_2^{-1} $. Then as $n \to\infty$, $\mathrm{ FDP}_\Phi, \mathrm{
FDP}_{\mathrm{ B}} \rightarrow^{ P} \alpha\pi_0$ and $\mathrm{ FDR}_\Phi,
\mathrm{
FDR}_{\mathrm{ B}} \to\alpha\pi_0$.
\end{theorem}

Attractive properties of the bootstrap for multiple-hypothesis testing
were first noted by \citet{Hall1990} in the case of the mean rather than
its Studentized counterpart. Now it has been rigorously proved that
bootstrap methods are particularly effective in relieving skewness in
the extreme tails which leads to second-order accuracy [\citet
{FanHallYao2007,DelaigleHallJin2011}]. It is interesting and challenging
to investigate whether these advantages of the bootstrap can be
inherited by multiple $U$-testing in either the standardized or the
Studentized case.

\section{Numerical study} \label{sim.sec}

In this section, we present numerical investigations for various
calibration methods described in Section~\ref{HT.sec} when they are
applied to two-sample large-scale multiple testing problems. We refer
to the simulation for two-sample $t$-test and Studentized Mann--Whitney
test as $\mathrm{Sim}_1$ and $\mathrm{Sim}_2$, respectively. Assume that we observe two
groups of $m$-dimensional gene expression data $\{\bX_i\}_{i=1}^{n_1}$
and $\{\bY_j\}_{j=1}^{n_2}$, where $\bX_{1},\ldots,\bX_{n_1}$ and
$\bY
_1,\ldots,\bY_{n_2}$ are independent random samples drawn from the
distributions of $\bX$ and $\bY$, respectively.

For $\mathrm{Sim}_1$, let $\bX$ and $\bY$ be such that
%
\begin{equation}
\mathbf{X} = \bolds{\mu}_1+ \bigl\{ \bolds{\varepsilon}_1-
\e(\bolds{\varepsilon}_1)\bigr\} \quad\mbox {and}\quad \mathbf{Y}= \bolds{
\mu}_2+ \bigl\{ \bolds{\varepsilon} _2 - \e(\bolds{
\varepsilon}_2 )\bigr\},
\end{equation}
where $\bolds{\varepsilon}_1=(\varepsilon_{1,1},\ldots
,\varepsilon_{1,m})^{\mathrm{T}}$ and $\bolds{\varepsilon
}_2=(\varepsilon
_{2,1},\ldots
,\varepsilon_{2,m})^{\mathrm{T}}$ are two sets of i.i.d. random
variables. The
i.i.d. components of noise vectors $\bolds{\varepsilon}_1$ and
$\bolds{\varepsilon}_2$ follow two
types of distributions: (i) the exponential distribution $\operatorname
{Exp}(\lambda)$ with density function $\lambda^{-1}e^{-x/\lambda}$;
(ii) Student $t$-distribution $t(k)$ with $k$ degrees of freedom. The
exponential distribution has nonzero skewness, while the
$t$-distribution is symmetric and heavy-tailed. For each type of error
distribution, both cases of homogeneity and heteroscedasticity were
considered. Detailed settings for the error distributions are specified
in Table~\ref{table:a1}.
%
\begin{table}[b]
\caption{Distribution settings in $\mathrm{Sim}_1$}\label{table:a1}
\begin{tabular*}{\textwidth}{@{\extracolsep{\fill}}lcc@{}}
\hline
& \textbf{Homogeneous case} & \textbf{Heteroscedastic case} \\
\hline
Exponential distributions & $\varepsilon_{1,k}\sim\operatorname{Exp}(2)$ &
$\varepsilon_{1,k}\sim\operatorname{Exp}(2)$ \\
& $\varepsilon_{2,k}\sim\operatorname{Exp}(2)$ & $\varepsilon_{2,k}\sim
\operatorname{Exp}(1)$ \\[3pt]
Student $t$-distributions & $\varepsilon_{1,k}\sim t(4)$ &
$\varepsilon
_{1,k}\sim t(4)$\\
& $\varepsilon_{2,k}\sim t(4)$ & $\varepsilon_{2,k}\sim t(3)$\\
\hline
\end{tabular*}
\end{table}

For $\mathrm{Sim}_2$, we assume that $\bX$ and $\bY$ satisfy
%
\begin{equation}
\mathbf{X} = \bolds{\mu}_1+ \bolds{\varepsilon}_1\quad
\mbox{and}\quad \mathbf{Y}= \bolds{\mu}_2+ \bolds{\varepsilon}_2
,
\end{equation}
where
$\bolds{\varepsilon}_1=(\varepsilon_{1,1},\ldots,\varepsilon
_{1,m})^{\mathrm{T}}$ and $\bolds{\varepsilon}
_2=(\varepsilon_{2,1},\ldots,\varepsilon_{2,m})^{\mathrm{T}}$ are
two sets of
i.i.d. random variables. We consider several distributions for the
error terms $\varepsilon_{1,k}$ and $\varepsilon_{2,k}$: standard
normal distribution $N(0,1)$, $t$-distribution $t(k)$, uniform
distribution $U(a,b)$ and Beta distribution $\operatorname{Beta}(a,b)$.
Table~\ref{table:a2} reports four settings of $(\varepsilon_{1,k},\varepsilon
_{2,k})$ used
in our simulation. In either setting, we know $\p(\varepsilon
_{1,k}\leq
\varepsilon_{2,k})=1/2$ holds. Hence, the power against the null
hypothesis $H_0^k:\p(X_{k}\leq Y_k)=1/2$ will generate from the
magnitude of the difference between the $k$th components of $\bolds
{\mu}_1$
and $\bolds{\mu}_2$.
%
\begin{table}[t]
\caption{Distribution settings in $\mathrm{Sim}_2$}\label{table:a2}
\begin{tabular*}{\textwidth}{@{\extracolsep{\fill}}lcc@{}}
\hline
& \textbf{Identical distributions} & \textbf{Nonidentical distributions} \\
\hline
Case 1 & $\varepsilon_{1,k}\sim N(0,1)$ & $\varepsilon_{1,k}\sim
N(0,1)$ \\
& $\varepsilon_{2,k}\sim N(0,1)$ & $\varepsilon_{2,k}\sim t(3)$ \\[3pt]
Case 2 & $\varepsilon_{1,k}\sim U(0,1)$ & $\varepsilon_{1,k}\sim
U(0,1)$\\
& $\varepsilon_{2,k}\sim U(0,1)$ & $\varepsilon_{2,k}\sim\operatorname
{Beta}(10,10)$\\
\hline
\end{tabular*}
\end{table}

\begin{figure}

\includegraphics{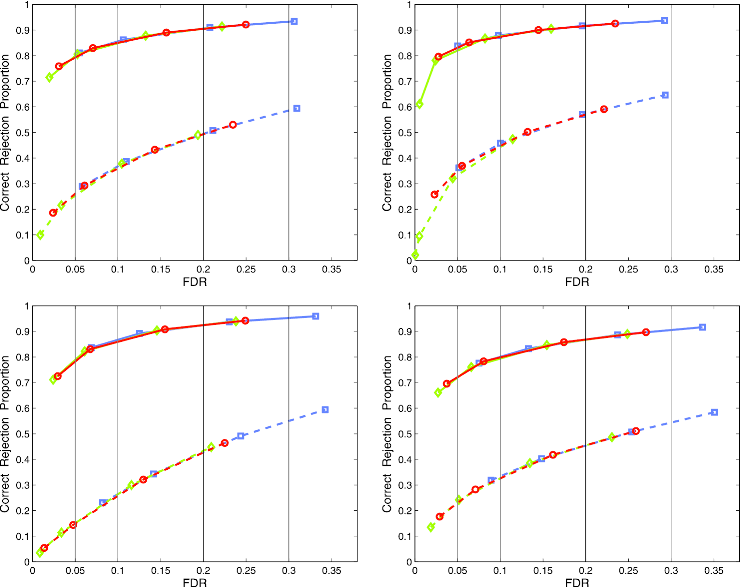}

\caption{Performance comparison of B--H procedures based on three
calibration methods in $\mathrm{Sim}_1$ with $(n_1,n_2)=(50,30)$ and $m=1000$.
The first and second rows show the results when the components of noise
vectors $\bolds{\varepsilon}_1$ and $\bolds{\varepsilon}_2$ follow
$t$-distributions and
exponential distributions, respectively; left and right panels show the
results for homogeneous and heteroscedastic cases, respectively;
horizontal and vertical axes depict empirical false discovery rate and
empirical correct rejection proportion, respectively; and the
prescribed levels $\alpha=0.05, 0.1, 0.2$ and $0.3$ are indicated by
unbroken horizontal black lines. In each panel, dashed lines and
unbroken lines represent the results for the discrepancy parameter
$c=1$ and $1.5$, respectively, and different colors express different
methods employed to calculate $p$-values in the B--H procedure, where
blue line, green line and red line correspond to the procedures based
on normal, conventional and regularized bootstrap calibrations,
respectively.}\label{fig1}
\end{figure}

In both $\mathrm{Sim}_1$ and $\mathrm{Sim}_2$, we set $\bolds{\mu}_1=\bzero$, and
assume that
the first $m_1=\lfloor1.6m^{1/2}\rfloor$ components of $\bolds{\mu
}_2$ are
equal to $c\{(\sigma_1^2 n_1^{-1} + \sigma_2^2 n_2^{-1} )\log m \}
^{1/2}$ and the rest are zero. Here, $\sigma_1^2$ and $\sigma_2^2$
denote the variance of $\varepsilon_{1,k}$ and $\varepsilon_{2,k}$, and
$c$ is a parameter employed to characterize the location discrepancy
between the distributions of $\bX$ and~$\bY$. The sample size $(n_1,
n_2)$ was set to be $(50,30)$ and $(100,60)$, and the discrepancy
parameter $c$ took values in $\{1, 1.5 \}$. The significance level
$\alpha$ in the B--H procedure was specified as $0.05,0.1,0.2$ and
$0.3$, and the dimension $m$ was set to be $1000$ and $2000$. In
$\mathrm{Sim}_1$, we compared three different methods to calculate the
$p$-values in the B--H procedure: normal calibration given in Section~\ref{se:norm}, bootstrap calibration and regularized bootstrap
calibration proposed in
Section~\ref{reg.bootstrap}. For regularized bootstrap calibration, we
used a cross-validation approach as in Section~3 of \citet{LiuShao2014}
to choose regularized parameters $\la_{1k}$ and $\la_{2k}$. In $\mathrm{Sim}_2$,
we compared the performance of normal calibration and bootstrap
calibration proposed in Section~\ref{app.MW}. For each compared method,
we evaluated its performance via two indices: the empirical FDR and the
proportion among the true alternative hypotheses was rejected. We call
the latter correct rejection proportion. If the empirical FDR is low,
the proposed procedure has good FDR control; if the correct rejection
proportion is high, the proposed procedure has fairly good performance
in identifying the true signals. For ease of exposition, we only report
the simulation results for $(n_1,n_2)=(50,30)$ and $m=1000$ in Figures
\ref{fig1} and \ref{fig2}.
The results for $(n_1,n_2)=(100,60)$ and $m=2000$ are similar,
which can be found in the supplementary material [\citet{supp}]. Each curve corresponds
to the performance of a certain method and the line types are specified
in the caption below. The horizontal ordinates of the four points on
each curve depict the empirical FDR of the specified method when the
pre-specified level $\alpha$ in the B--H procedure was taken to be
$0.05, 0.1, 0.2$ and $0.3$, respectively, and the vertical ordinates
indicate the corresponding empirical correct rejection proportion. We
say that a method has good FDR control if the horizontal ordinates of
the four points on its performance curve are less than the prescribed
$\alpha$ levels.

\begin{figure}

\includegraphics{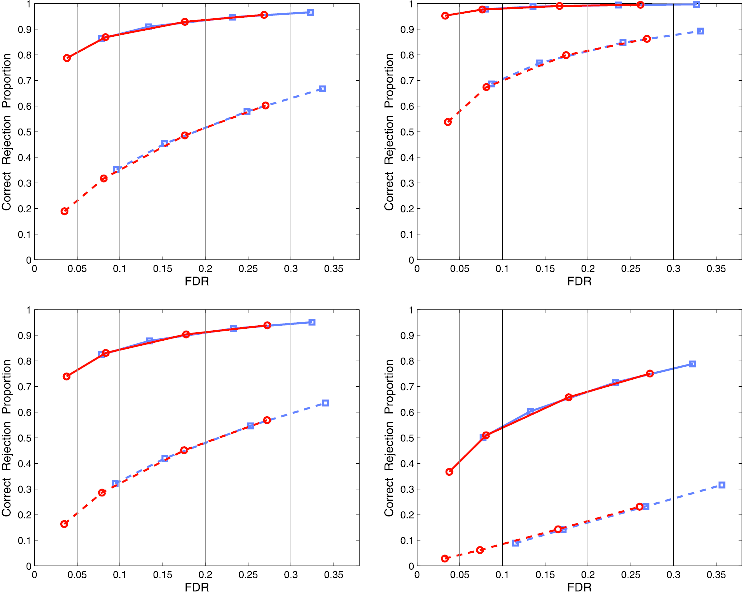}

\caption{Performance comparison of B--H procedures based on two
different calibration methods in $\mathrm{Sim}_2$ with $(n_1,n_2)=(50,30)$ and
$m=1000$. The first and second rows show the results when the
components of noise vectors $\bolds{\varepsilon}_1$ and $\bolds
{\varepsilon}_2$ follow the
distributions specified in cases 1 and 2 of Table~\protect\ref{table:a2}, respectively;
left and right panels show the results for the cases of identical
distributions and nonidentical distributions, respectively; horizontal
and vertical axes depict empirical false discovery rate and empirical
correct rejection proportion, respectively; and the prescribed levels
$\alpha=0.05, 0.1, 0.2$ and $0.3$ are indicated by unbroken horizontal
black lines. In each panel, dashed lines and unbroken lines represent
the results for the discrepancy parameter $c = 1$ and $1.5$,
respectively, and different colors express different methods employed
to calculate $p$-values in the B--H procedure, where blue line and red
line correspond to the procedures based on normal and bootstrap
calibrations, respectively.}\label{fig2}
\end{figure}

In general, as shown in Figures \ref{fig1} and \ref{fig2}, the B--H
procedure based on
(regularized) bootstrap calibration has better FDR control than that
based on normal calibration. In $\mathrm{Sim}_1$ where the errors are symmetric
(e.g., $\varepsilon_{1,k}$ and $\varepsilon_{2,k}$ follow the Student
$t$-distributions), the panels in the first row of Figure~\ref{fig1} show that
the B--H procedures using all the three calibration methods are able to
control or approximately control the FDR at given levels, while the
procedures based on bootstrap and regularized bootstrap calibrations
outperform that based on normal calibration in controlling the FDR.
When the errors are asymmetric in $\mathrm{Sim}_1$, the performances of the
three B--H procedures are different from those in the symmetric cases.
From the second row of Figure~\ref{fig1}, we see that the B--H procedure based on
normal calibration is distorted in controlling the FDR while the
procedure based on (regularized) bootstrap calibration is still able to
control the FDR at given levels. This phenomenon is further evidenced
by Figure~\ref{fig2} for $\mathrm{Sim}_2$. Comparing the B--H procedures based on
conventional and regularized bootstrap calibrations, we find that the
former approach is uniformly more conservative than the latter in
controlling the FDR. In other words, the B--H procedure based on
regularized bootstrap can identify more true alternative hypotheses
than that using conventional bootstrap calibration. This phenomenon is
also revealed in the heteroscedastic case. As the discrepancy parameter
$c$ gets larger so that the signal is stronger, the correct rejection
proportion of the B--H procedures based on all the three calibrations
increase and the empirical FDR is closer to the prescribed level.

\section{Discussion}\label{se:discu}

In this paper, we established Cram\'er-type moderate deviations for
two-sample Studentized $U$-statistics of arbitrary order in a general
framework where the kernel is not necessarily bounded. Two-sample
$U$-statistics, typified by the two-sample Mann--Whitney test
statistic, have been widely used in a broad range of scientific
research. Many of these applications rely on a misunderstanding of what
is being tested and the implicit underlying assumptions, that were not
explicitly considered until relatively recently by \citet
{ChungRomano2011}. More importantly, they provided evidence for the
advantage of using the Studentized statistics both theoretically and
empirically.

Unlike the conventional (one- and two-sample) $U$-statistics, the
asymptotic behavior of their Studentized counterparts has barely been
studied in the literature, particularly in the two-sample case.
Recently, \citet{ShaoZhou2012} proved a Cram\'er-type moderate deviation
theorem for general Studentized nonlinear statistics, which leads to a
sharp moderate deviation result for Studentized one-sample
$U$-statistics. However, extension from one-sample to two-sample in the
Studentized case is totally nonstraightforward, and requires a more
delicate analysis on the Studentizing quantities. Further, for the
two-sample $t$-statistic, we proved moderate deviation with
second-order accuracy under a finite $4$th moment condition (see
Theorem~\ref{thmm.t-stat}), which is of independent interest. In
contrast to the one-sample case, the two-sample $t$-statistic cannot
be reduced to a self-normalized sum of independent random variables,
and thus the existing results on self-normalized ratios [\citet
{JingShaoWang2003}, \citeauthor{Wang2005} (\citeyear{Wang2005,Wang2011})] cannot be directly applied.
Instead, we modify Theorem~2.1 in \citet{ShaoZhou2012} to obtain a more
precise expansion that can be used to derive a refined result for the
two-sample $t$-statistic.

Finally, we show that the obtained moderate deviation theorems provide
theoretical guarantees for the validity, including robustness and
accuracy, of normal, conventional bootstrap and regularized bootstrap
calibration methods in multiple testing with FDR/FDP control. The
dependence case is also covered. These results represent a useful
complement to those obtained by \citet{FanHallYao2007}, \citet
{DelaigleHallJin2011} and \citet{LiuShao2014} in the one-sample case.

\section*{Acknowledgements}
The authors would like to thank Peter Hall and Aurore Delaigle for
helpful discussions and encouragement. The authors sincerely thank the
Editor, Associate Editor and three referees for their very constructive
suggestions and comments that led to substantial improvement of the paper.

%
\begin{supplement}[id=suppA]
\stitle{Supplement to ``Cram\'er-type moderate deviations for
Studentized two-sample $U$-statistics with applications''}
\slink[doi]{10.1214/15-AOS1375SUPP} 
\sdatatype{.pdf}
\sfilename{aos1375\_supp.pdf}
\sdescription{This supplemental material contains proofs for all the
theoretical results in the main text, including Theorems \ref{thmm.1},
\ref{thmm.t-stat}, \ref{app.thmm.1}
and \ref{FDR.dependence}, and additional numerical results.}
\end{supplement}

%




\printaddresses
\end{document}